\documentclass{article}

\usepackage{amsmath,amsthm,amssymb,amscd}
\usepackage{amsfonts}
\usepackage{rotating}
\usepackage{euscript}
\usepackage{pst-node}
\usepackage{epsfig,verbatim}
\usepackage{pictexwd,dcpic}

\def\be{\begin{equation}}
\def\ee{\end{equation}}

\def\C{{\mathbb C}} 
\def\f{\EuScript}
 
\def\P{{\mathbb P}}
\def\Z{{\mathbb Z}}

\def\e{\eqref}
\def\phi{{\varphi}}
 
\def\tt{\widetilde}
\def\deg{{\rm deg\,}}

\def\cos{{\rm cos\,}} 
 
\def\GCD{{\rm GCD }}
\def\LCM{{\rm LCM }}

\def\bp{\begin{proposition}}
\def\ep{\end{proposition}}

\def\bt{\begin{theorem}}
\def\et{\end{theorem}}
\def\br{\begin{remark}}
\def\er{\end{remark}}
\def\be{\begin{equation}}
\def\bee{\begin{equation*}}
\def\l{\label}

\def\e{\eqref}

\def\ee{\end{equation}}
\def\eee{\end{equation*}}
\def\bl{\begin{lemma}}
\def\el{\end{lemma}}
\def\bc{\begin{corollary}}
\def\ec{\end{corollary}}
\def\pr{\noindent{\it Proof. }}

\def\bd{\begin{definition}}
\def\ed{\end{definition}}
\def\t{\widetilde}

\def\bco{\begin{conjecture}}
\def\eco{\end{conjecture}}
\def\bpr{\begin{problem}}
\def\epr{\end{problem}}

\newtheorem{theorem}{Theorem}[section]
\newtheorem{lemma}{Lemma}[section]
\newtheorem{definition}{Definition}[section]
\newtheorem{corollary}{Corollary}[section]
\newtheorem{proposition}{Proposition}[section]

\newtheorem{conjecture}{Conjecture}[section]
\newtheorem{problem}{Problem}[section]
\theoremstyle{definition}
\newtheorem{remark}[theorem]{Remark}

\begin{document}
\title{On semiconjugate rational functions}
\author{F. Pakovich}
\date{\ }
\maketitle

\begin{abstract} 
We investigate semiconjugate rational functions, that is rational functions $A,$ $B$ related by the functional equation 
$A\circ X=X\circ B$, where $X$ is a rational function. We show that if $A$ and $B$ is a pair of such functions, then either $A$ can be obtained from $B$ by a certain iterative process, or $A$ and $B$ can be described in terms of orbifolds of non-negative 
Euler characteristic 
on the Riemann sphere.

\end{abstract}

\section{Introduction} 
Let $A$, $B$ be two rational functions of degree at least two on the Riemann sphere. 
The function $B$ is said to be semiconjugate to the function $A$ if there exists a non-constant rational function $X$
such that the equality 
\be \l{i1}
A\circ X=X\circ B
\ee
holds. If $X$ is invertible the functions $A$ and $B$ are called conjugate.  The semiconjugacy relation plays an important role in the study of complex dynamical systems, and in this article we study triples  
$A,B,X$ which satisfy relation \eqref{i1}.

Since condition \eqref{i1} is not symmetric with respect to $A$ and $B$, the semiconjugacy is not an equivalency relation. However, it is easy to see that if  $B$ is semiconjugate to $A$, and $A$ is semiconjugate to $C$, then 
$B$ is semiconjugate to $C$. Therefore, the semiconjugacy is a preorder on the set of rational functions. 
Having this in mind, we will use the notation $A\leq B$ for rational functions  $A$ and $B$ satisfying equality \eqref{i1} for some rational function $X.$

The problem of describing of semiconjugate rational functions can be considered as a generalization  of 
the classical problem of 
describing of commuting rational functions, that is of rational functions  $A$ and $X$ satisfying the functional equation 
\be \l{comm} A\circ X=X\circ A.\ee
The last problem was considered  in the early twenties of the past century in the papers of Fatou, Julia, and Ritt 
\cite{f}, \cite{j}, \cite{r}. 
In all these papers it was  
assumed that the considered commuting functions $A$ and $X$ have no 
iterate in common, that is 
\be \l{ins}  A^{\circ n}\neq X^{\circ m}\ee for all $n,m\in \mathbb N.$
In particular, this assumption rules out ``trivial'' solutions of the form $A=R^{\circ m},$ $X=R^{\circ n},$
where $R$ is an arbitrary rational function.
Fatou and Julia
used dynamical methods requiring an additional assumption that the
Julia set of $A$ or $X$ does not coincide with the whole complex plane, while
Ritt
used a method of algebraic-topological character free of any assumptions about the Julia set.
Briefly, the Ritt theorem states that if rational functions $A$ and $X$ commute and
no iterate of $A$ is equal to an iterate of $X$,
then, up to a conjugacy, $A$ and $X$ are either powers, or  Chebyshev polynomials, or Latt\`es functions.  

Notice that both equations \eqref{i1} and \eqref{comm} are particular cases of the more general functional equation 
\be \l{rii} A\circ X=Y\circ B\ee
investigated for the first time by Ritt in the paper \cite{r2}. In this paper Ritt laid the foundation of the decomposition theory of rational functions, and constructed a comprehensive decomposition theory of polynomials.
In particular, Ritt described solutions of \eqref{rii} in the case where  $A,B,X,Y$ are polynomials, and the results of \cite{r2} can be applied to equation \eqref{i1}  
in the polynomial case (see \cite{i}).
The Ritt theory may be extended to a decomposition theory of Laurent polynomials  (\cite{pak}). However, more general results in this direction are not known.

A proof of the Ritt theorem, based on modern dynamical 
methods,
was given by Eremenko \cite{e2} who pointed out that the Ritt result can be formulated in a natural way using  the concept of orbifold introduced by Thurston 
\cite{t}. Recall that a Riemann surface orbifold is a pair $\f O=\f(R,\nu)$ consisting of a Riemann surface $R$ and a ramification function $\nu:R\rightarrow \mathbb N$ which takes the value $\nu(z)=1$ except at isolated set of points. The Euler characteristic of an orbifold $\f O=(R,\nu)$ is defined by the formula  
\be\l{frh} \chi(\f O)=\chi(R)+\sum_{z\in R}\left(\frac{1}{\nu(z)}-1\right),\ee where $\chi(R)$ is the 
Euler characteristic of $R.$
If $R_1$, $R_2$ are Riemann surfaces provided with ramification functions $\nu_1,$ $\nu_2$, and 
$f:\, R_1\rightarrow R_2$ is a holomorphic branched covering map, then $f$
is said to be a covering map $f:\,  \f O_1\rightarrow \f O_2$
between orbifolds
$\f O_1=(\f R_1,\nu_1)$ and $\f O_2=\f(R_2,\nu_2)$
if for any $z\in R_1$ the equality 
\be \l{us+} \nu_{2}(f(z))=\nu_{1}(z)\deg_zf\ee holds, where $\deg_zf$ is the local degree of $f$ at the point $z$. 
If such a map has a finite degree $d$, then the Riemann-Hurwitz 
formula implies that 
\be \l{rhor+} \chi(\f O_1)=d \chi(\f O_2). \ee

In the above terms the Ritt theorem may be formulated as follows (\cite{e2}):  if rational functions $A$ and $X$ commute and
no iterate of $A$ is equal to an iterate of $X$,
then there exists an orbifold $\f O=(R,\nu)$ of {\it zero} Euler characteristic  with
$R$ equal  to $\C\setminus\{0\}$, $\C$, or $\C\P^1$ such that 
the commutative diagram 
$$
\begin{CD}
{\f O} @>A>> {\f O}\\
@VV X V @VV X V\\ 
{\f O} @>A >> {\f O}\ 
\end{CD}
$$
consists of covering  maps between orbifolds. In this description power functions $z^{\pm n}$ correspond to covering maps preserving the orbifold 
$\f O=(R,\nu)$, where $R=\C\setminus\{0\}$ and $\nu\equiv 1.$ The Chebyshev polynomials $\pm T_n$, defined by the formula $T_n(\cos \phi)=\cos(n \phi)$, 
correspond to covering maps preserving the orbifold 
$\f O=(R,\nu)$, where $R=\C$ and $\nu(-1)=\nu(1)=2.$
Finally, Latt\`es functions correspond to covering maps preserving an orbifold of zero 
characteristic $\f O$ with $R=\C\P^1$. In the last case formula \eqref{frh} implies that the collection of ramification indices of $\f O$ is either $(2,2,2,2)$, or one of the following triples $(3,3,3)$, $(2,4,4)$, $(2,3,6)$. 
Notice that the Ritt theorem provides no information about functions (commuting or not) that do share an iterate, and  
a description of such functions 
is known only in the polynomial case (\cite{r3}, \cite{r}).  Thus, in a certain sense the classification of commuting rational functions is not yet completed.

In comparison with equation \eqref{comm} equation \eqref{i1} has many more solutions.
Indeed, take arbitrary rational functions $U_1,V_1$ and set 
$$B=V_1\circ U_1, \ \ \ A=U_1\circ V_1.$$ 
Then the equality 
\be \l{ii} (U_1\circ V_1)\circ U_1= U_1\circ (V_1\circ U_1)\ee
implies that $A\leq B$. Similarly, 
$B\leq A$. Moreover, if  now
$U_2,V_2$ are rational functions such that the equality 
\be \l{ko}  U_{1}\circ V_{1}=V_{2}\circ U_{2}\ee holds, 
then the function $H=U_2\circ V_2$ satisfies $A\leq H$ and
$H\leq A$, implying that $B\leq H$ and
$H\leq B.$

This motivates the following definition of an equivalence relation  on the set of  rational functions:
$B\sim A$ if there exist rational functions $U_i,V_i,$ $1\leq i \leq n,$ such that $B=V_1\circ U_1,$ 
\be \l{xori0} U_i\circ V_i=V_{i+1}\circ U_{i+1}, \ \ \ 1\leq i \leq n-1,\ee and 
$A=U_n\circ V_n.$ 
Notice that since for any rational function $W$ of degree one the equality
$$B=(B\circ W)\circ W^{-1}$$ holds, each equivalence class of $\sim$ is a union of conjugacy classes. Thus, since $B\sim A$ implies that $A\circ X=X\circ B$ and $B\circ Y=Y\circ  A$ 
for 
\be \l{xori} X=U_n\circ U_{n-1} \circ \dots \circ U_1  \ee
and
$$Y=V_1\circ V_{2} \circ \dots \circ V_n,$$ the equivalence  relation $\sim$ can be considered as a weaker form of the classical conjugacy relation whose classes consist of  functions   having  ``similar'' although not ``identical'' dynamics.

Roughly speaking, our main result states that unless $B\sim A$ the relation $A\leq B$ implies very strong restrictions
on $A$ and $B,$ which can  
be described in terms of orbifolds of {\it non-negative} Euler characteristic on the Riemann sphere. 
Namely, similarly to Latt\`es functions, such $A$ and $B$ can be characterized as maps ``preserving'' some orbifold on the Riemann sphere.
However, since \eqref{rhor+} implies that $\chi(\f O)= 0$ for any self-covering map $f:\, \f O\rightarrow \f O$, we take as a basis the following weakened modification of the notion of covering. 

A rational function $f$ is called {\it a  holomorphic map} $f:\, \f O_1\rightarrow \f O_2$
between orbifolds
$\f O_1$ and $\f O_2$ defined on $\C\P^1$ if for any $z\in \C\P^1$ instead of equality \eqref{us+} 
a weaker condition  
\be \l{vbn} \nu_{2}(f(z)) \mid \nu_{1}(z)\deg_zf\ee holds. 
For fixed $\nu_2(z)$ a minimal possible value for $\nu_1(z)$ such that \eqref{vbn} holds is defined by 
the equality 
\be \l{rys+} \nu_{2}(f(z))=\nu_{1}(z)\GCD(\deg_zf, \nu_{2}(f(z)).\ee 
If \eqref{rys+} is satisfied for all $z\in \C\P^1$ we 
say that $f$ is {\it a minimal holomorphic  map} 
between orbifolds. Notice that   any covering map between orbifolds is necessarily a minimal holomorphic  map.
The importance of the notion of minimal holomorphic map is explained by the fact that for any $f$ and   $\f O$ 
there exists a uniquely defined orbifold $f^*\f O$ such that 
$f:\, f^*\f O\rightarrow \f O$ is a  minimal holomorphic map, and   the equality  $(g\circ f)^*\f O= f^*(g^*\f O)$ holds for any $f$, $g$, and   $\f O$.

If $f:\, \f O_1\rightarrow \f O_2$ is a  holomorphic map between orbifolds, then instead of equality \eqref{rhor+} the inequality 
\be \l{iioopp+} \chi(\f O_1)\leq \chi(\f O_2)\,\deg f \ee holds, and the equality 
is attained if and only if $f:\, \f O_1\rightarrow \f O_2$ is a covering map. In particular, this implies that if 
$f:\, \f O\rightarrow \f O$
is a minimal holomorphic {\it self-map} and $\deg f>1$, then $\chi(\f O)\geq 0$. Moreover, if $\chi(\f O)= 0$, then $f$ in fact is  a covering self-map and hence a Latt\`es function. Thus, the class of {\it minimal holomorphic self-maps} between orbifolds is a natural extension of the class
of Latt\`es functions. Notice that although we assume that all considered orbifolds are defined on $\C\P^1$, the functions $\pm z^n$ and $\pm T_n$ still belong to  this new class. Indeed, it is easy to see 
that $z^{\pm n}:\f O\rightarrow \f O$ is a minimal holomorphic map between orbifolds for any $\f O$ defined by the conditions 
$$\nu(0)=\nu(\infty)=m,\ \ \ \GCD(n,m)=1,$$ while  $\pm T_{n}:\f O\rightarrow \f O$ is a minimal holomorphic map between orbifolds for any $\f O$ defined by the conditions 
$$\nu(-1)=\nu(1)=2, \ \ \ \nu(\infty)=m, \ \ \ \GCD(n,m)=1.$$

In the above notation our main result may be formulated as follows.

\bt \l{main} Let $A,B$ and $X$ be rational functions of degree at least two  such that  $A\circ X=X\circ B$.
Then either $B\sim A$ and $X$ satisfies \eqref{xori} for some chain of transformations \eqref{xori0}, or 
there exist orbifolds $\f O_1$, $\f O_2$ of non-negative Euler characteristic on the Riemann sphere such that 
the commutative diagram 
$$
\begin{CD}
{\f O_1} @>B>> {\f O_1}\\
@VV X V @VV X V\\ 
{\f O_2} @>A >> {\f O_2}\ 
\end{CD}
$$
consists of  minimal holomorphic  maps between orbifolds. Furthermore, either $\chi(\f O_1)= \chi(\f O_2)=0$
and $A$, $B$ are Latt\`es functions, or $0<\chi(\f O_2)< \chi(\f O_1)$. In the last case, the possible  
collections of ramification indices of $\f O_1$ and $\f O_2$ are the following: $(n,n)$ or $(2,2,n)$ for some $n\geq 2,$ or one of the triples $(2,3,3)$, $(2,3,4),$ $(2,3,5).$ 
In addition, $\f O_1$ may be a non-ramified sphere. 

\et

For example, for the solution  
$$T_n\circ \frac{1}{2}\left(z^m+\frac{1}{z^m}\right)=\frac{1}{2}\left(z^m+\frac{1}{z^m}\right)\circ z^n$$ 
of \eqref{i1} the orbifold $\f O_2$ is defined by the conditions 
\be \l{barsi} \nu(-1)=\nu(1)=2, \ \ \ \nu(\infty)=m/d,\ee
where $d=\GCD(n,m)$, while $\f O_1$ is non-ramified. 
On the other hand, for the solution  
$$T_n\circ T_m= T_m\circ T_n$$ the orbifold $\f O_2$ is still defined by \eqref{barsi} but 
$\f O_1$ is defined  by  the condition
$$\nu(-1)=\nu(1)=2.$$ Notice that in the last example $A=B=T_n$ so $A\sim B$. However, the corresponding $X=T_m$ in general cannot be obtained 
by formula \eqref{xori} so the second part of the theorem is applied.

If  $B\not\sim A$, then the  
conditions imposed by Theorem \ref{main} on $A,B,$ and $X$ are quite strong and provide a reasonably good description of solutions of \eqref{i1}. Different implications of these conditions are discussed in the body of the paper. 
In contrast, the paper contains essentially no information about the relation $A\sim B$.
Notice however that this relation also is rather restrictive. Indeed,  
it follows easily from the definition, that any two equivalent rational functions are {\it isospectral} in the following sense: 
for each $n \geq 1$ the unordered lists of multipliers at fixed points of the iterates $A^{\circ n}$ and $B^{\circ n}$ are the same.
By the result of McMullen (\cite{Mc}), this implies in particular that unless $A$ is a flexible Latt\`es function its equivalence class contains at most a {\it finite} number of conjugacy classes.

The paper is organized as follows. 
In the second section we review  general properties of the functional equation 
\be \l{l} f\circ p=g\circ q\ee on compact Riemann surfaces, basing on the fiber product 
approach. 
In the third section  we provide main definitions and results concerning orbifolds on 
Riemann surfaces and holomorphic maps between orbifolds.

In the fourth section we introduce the concept of a minimal holomorphic map between orbifolds, and 
study properties of such maps with respect to functional decompositions.
In particular, we 
relate minimal holomorphic maps with functional equation \e{l}.
In the fifth section we describe an approach to the classification of minimal holomorphic self-maps between orbifolds of positive Euler characteristic on the Riemann sphere. 
In particular, we relate such maps 
with rational functions equivariant with respect to finite 
subgroups of $Aut(\C\P^1)$.
Finally, in the sixth section  we prove  Theorem \ref{main} 
and provide a number of  examples.

\section{Equation $f\circ p=g\circ q$ and fiber products}
In this section we recall, mostly without proofs, some general results related to the functional equation
\be \l{m} h=f\circ p=g\circ q, \ee 
where 
$h:\, R\rightarrow \C\P^1,$ $p:\, R\rightarrow C_1$, $f:\, C_1\rightarrow \C\P^1,$ $q:\, R\rightarrow C_2$, $g:\, C_2\rightarrow \C\P^1$ 
are holomorphic functions on compact Riemann surfaces.
For more details we refer the reader to \cite{pak}, Section 2 and 3.

Let $h:\, R\rightarrow \C\P^1$ be a holomorphic function on a compact Riemann surface $R$, and 
$\f C(h)=\{z_1, z_2, \dots , z_r\}$ the set of critical values of $h$. Fix a point \linebreak
$z_0\in\C\P^1\setminus \f C(h)$ and some loops $\gamma_i$ around $z_i,$ $1\leq i \leq r,$ such that 
$\gamma_1\gamma_2...\gamma_r=1$ in $\pi_1(\C\P^1\setminus \f C(h),z_0)$. 
Denote by 
$\delta_i,$ $1\leq i \leq r,$ a permutation of points of 
$h^{-1}\{z_0\}$ induced by the lifting of $\gamma_i,$ $1\leq i \leq r,$ by $h$.
The permutation group $G_h$ generated by $\delta_i,$ $1\leq i \leq r,$  is called the monodromy group of $h$. 
Clearly, 
\be \l{posl} \delta_1\delta_2...\delta_r=1\ee 
in $G_h$.

Let $G$ be a group which acts transitively on a finite set $X$. 
Recall that a subset $B$ of $X$ is called a block of  $G$, if for each $g\in G$ either $g(B)=B$ or $g(B)\cap B=\emptyset.$  
Clearly, if $B$ is a block, then $\f B=\{\sigma(B), \sigma\in G\}$ is a partition of $X$, which is called an imprimitvity system of $G$.  

The monodromy group $G_h$ is related to compositional properties of the function 
$h$ as follows.
If  $h$
can be decomposed into a composition $h=f\circ p$ of holomorphic 
functions $p:\, R\rightarrow C_1$ and $f:\, C_1\rightarrow \C\P^1$, where $\deg f=d,$ 
then $G_{h}$ has an imprimitivity system consisting of 
$d$ blocks $\f A_i=p^{-1}\{t_i\},$ $1\leq i \leq d,$ where $\{t_1,t_2,\dots, t_{d}\} =f^{-1}\{z_0\}$, and the permutation group induced by the action of $G_h$ on these blocks is permutation isomorphic to the group $G_f.$
Furthermore, any imprimitivity system of $G_h$ arises from a decomposition of $h$,
and decompositions $h=f\circ p$ and $h=g\circ q$, where $q:\, R\rightarrow C_2$, 
$g:\, C_2\rightarrow \C\P^1,$ lead to the same
imprimitivity system if and only there exists an isomorphism $\mu:\, C_2 \rightarrow C_1$
such that $$ f=g \circ \mu^{-1}, \ \ \ p=\mu \circ q.$$

We will say 
that two holomorphic functions $p:\, R\rightarrow C_1$ and $q:\, R\rightarrow C_2$ have no non-trivial compositional common right factor, if 
the equalities 
\be \l{kaban} p= \tt p\circ  w, \ \ \ q= \tt q\circ w,\ee where $w:\, R \rightarrow \tt R$, $\tt p:\, \tt R\rightarrow C_1$, $\tt q:\, \tt R\rightarrow C_2$ are holomorphic functions, imply that $\deg w=1.$

\bt \l{es} 
For any two fixed holomorphic functions $f:\, C_1\rightarrow \C\P^1$ and $g:\, C_2\rightarrow \C\P^1$ 
there exist holomorphic functions $h_j:\, R_j\rightarrow \C\P^1$
(components of the fiber product of $f$ and $g$) and $p_j:\, R_j\rightarrow C_1$, $q_j:\, R_j\rightarrow C_2$
such that the following conditions are satisfied: 
\begin{itemize}
\item 
$h_j=f\circ p_j=g\circ q_j$, 
\item $\sum_{j}\deg h_j= \deg f\,\deg g,$ 
\item 
for any solution $h,p,q$ of \e{m} 
there exist an index $j$ and 
a holomorphic function $w:\, R\rightarrow R_j$ such that
$ h= h_j\circ w,$ $ p= p_j\circ  w,$ $q= q_j\circ w.$ 
\qed
\end{itemize} 
\et

Recall briefly the construction of the functions $h_j$. 
Let 
$S=\{z_1, z_2, \dots , z_r\}$ be the union of $\f C(f)$ and $\f C(g)$. As above, fix a point
$z_0$ from $\C\P^1\setminus S$ and small loops $\gamma_i$ around $z_i,$ $1\leq i \leq r,$ such that 
$\gamma_1\gamma_2...\gamma_r=1$ in $\pi_1(\C\P^1\setminus S,z_0)$. Set $n=\deg f,$ $m=\deg g$,  
and denote by  
$\alpha_i\in S_n$ (resp. $\beta_i\in S_m$) a permutation of points of 
$f^{-1}\{z_0\}$ (resp. $g^{-1}\{z_0\}$) induced by the lifting of $\gamma_i$, $1\leq i \leq r,$ by $f$ (resp. $g$). 
Clearly, the permutations $\alpha_i$ (resp. $\beta_i$), $1\leq i \leq r,$ generate the monodromy group of 
$f$ (resp. $g$) and \be \l{per} \alpha_1\alpha_2...\alpha_r=1,\ \ \ \ \beta_1\beta_2...\beta_r=1.\ee

Define now permutations $\delta_1, \delta_2, \dots, \delta_{r}\in S_{nm}$ as follows:
consider the set of $mn$ elements $c_{j_1,j_2},$ $1\leq j_1 \leq n,$
$1\leq j_2\leq m,$ and set $(c_{j_1,j_2})^{\delta_i}=c_{j_1^{\prime},j_2^{\prime}},$ where $$j_1^{\prime}
=j_1^{\alpha_i},\ \ \ \ j_2^{\prime}=j_2^{\beta_i}, \ \ \ \ 1\leq i \leq r.$$ 
It is convenient to consider $c_{j_1,j_2},$ $1\leq j_1 \leq n,$
$1\leq j_2\leq m,$ as elements of a $n\times m$ matrix $M$. Then the action of the 
permutation $\delta_i,$ $1\leq i \leq r,$ reduces 
to the permutation of rows of $M$ in accordance with the permutation $\alpha_i$ and 
the permutation of columns of $M$ in accordance with the permutation $\beta_i$.

In general, 
the permutation group generated by $\delta_i,$ $1\leq i \leq r,$
is not transitive on 
the set $c_{j_1,j_2},$ $1\leq j_1 \leq n,$
$1\leq j_2\leq m$. 
However, on each transitivity set $U_j$
the induced permutations $\delta_{i}(j),$ $1\leq i \leq r,$ 
satisfy the equality 
$$\delta_{1}(j)\delta_{2}(j)\dots \delta_{r}(j)=1.$$ 
By the Riemann existence theorem, 
this implies that there exist compact Riemann surfaces $R_j$ and holomorphic functions 
$h_j:\, R_j\rightarrow \C\P^1$ non-ramified outside $S$  
such that the permutations $\delta_{i}(j),$ 
$1\leq i \leq r,$ 
are induced by the lifting of $\gamma_i$ by $h_j.$ 
Moreover, it is easy to see by construction that 
the intersections of the transitivity set $U_j$ with the rows of $M$
form an imprimitivity system $\Omega_f(j)$ for the group generated by $\delta_{i}(j),$ $1\leq i \leq r,$  such that the permutations
of blocks of $\Omega_f(j)$ induced by 
$\delta_i(j),$ $1\leq i \leq r,$ coincide with $\alpha_i$. Similarly, 
the intersections of $U_j$ with the columns of $M$
form an imprimitivity system $\Omega_g(j)$ such that the permutations
of blocks of $\Omega_g(j)$ induced by 
$\delta_i(j),$ $1\leq i \leq r,$ coincide with $\beta_i.$ These imprimitivity systems correspond to decompositions $h_j=f\circ p_j=g\circ q_j$ for some
functions $p_j$ and $q_j$. 
\vskip 0.2cm

\bc \l{good2} Let $h,f,p,g,q$ be rational functions satisfying \e{m}. Then $p$ and $q$ have no non-trivial common compositional right factor if and only if 
for any $z\in \C\P^1$ the local degrees $\deg_zp$ and $\deg_zq$ are 
coprime. \ec 

\pr If $p$ and $q$ have no non-trivial common compositional right factor, then without loss of generality 
we can assume that $h=h_j,$ $p=p_j,$ $q=q_j$. 
Clearly,
\be \l{loi} \deg_zh_j=\deg_zp_j\,\deg_{p_j(z)}f=\deg_zq_j\,\deg_{q_j(z)}g.\ee On the other hand, the definition of 
the permutations $\delta_i,$ $1\leq i \leq r,$ yields that 
for any $z\in \C\P^1$ the equality  
\be \l{equ} \deg_zh_j=\LCM(\deg_{p_j(z)}f,\,\deg_{q_j(z)}g)\ee holds. 
It follows now from equalities \e{loi} and \e{equ} that $\deg_zp_j$ and $\deg_zq_j$ are 
coprime. 

In the other direction, if $w$ is a common compositional right factor of  $p$ and $q$, and  $z_0$ is any critical point of $w$, then the chain rule implies that $\deg_{z_0}p$ and $\deg_{z_0}q$ have a non-trivial common divisor. 
Therefore, since any rational function of degree greater than one has critical points, if 
$\deg_zp$ and $\deg_zq$ are 
coprime for any $z\in \C\P^1$,  $p$ and $q$ may not have a non-trivial common compositional right factor.
\qed

\vskip 0.2cm
Let $h,f,p,g,q$ be a solution of \e{m}, and  
$\f A$ and $\f B$ imprimitivity systems 
of the group $G_h$ corresponding to the decompositions $h=f\circ p$ and $h=g\circ q$ respectively. We say that 
the solution $h,f,p,g,q$ is {\it good} if any block of $\f A$ intersects with any block of $\f B$ and this intersection
consists of a unique element. This is equivalent to the requirement that the fiber product of $f$ and $g$ has a unique component
and $p$ and $q$ have no non-trivial common compositional right factor. 

Notice that a solution of \e{m}  in  rational functions $h,f,p,g,q$  is good if and only if the algebraic curve 
$$\f E(f,g):\, f(x)-g(y)=0$$ is irreducible and $\C(p,q)=\C(z)$. Indeed, 
irreducible components of $\f E(f,g)$  correspond to irreducible components of the
fiber product of $f$ and $g$  (see \cite{pak},  Proposition 2.4). On the other hand, by the L\"uroth theorem,
any subfield $K \subset \C(z)$, $K\neq \C$, has the form $K=\C(w)$ for some $w\in \C(z),$ implying that \eqref{kaban} holds for some rational functions $\t p,\t q, w$ with $\deg w>1$ if and only if $\C(p,q)\neq \C(z).$  

The construction of the fiber product 
implies easily the following statement.

\bl \l{good} A 
solution $h,f,p,g,q$ of \e{m} is good whenever 
any two of the following three conditions are satisfied:

\begin{itemize}
\item the fiber product of $f$ and $g$ has a unique component,
\item $p$ and $q$ have no non-trivial common compositional right factor,
\item $ \deg f=\deg q, \ \ \ \deg g=\deg p.$  \qed
\end{itemize}
\el

Finally, let us mention the following property 
of good solutions of \e{m}.

\bl \l{good1} Let $h,f,p,g,q$ be a good solution of \e{m},
$z_1$ a point from the set $g^{-1}\{z_0\}$, and 
$\sigma\in G_h$ a permutation which maps the set $q^{-1}\{z_1\}$ to itself. Then the permutation induced by $\sigma$ on $q^{-1}\{z_1\}$ has the same cyclic structure as 
the permutation induced by $\sigma$ on 
blocks $\f A_i=p^{-1}\{t_i\},$ $1\leq i \leq d,$ where $\{t_1,t_2,\dots, t_{d}\} =f^{-1}\{z_0\}.$
\el

\pr Since the set $q^{-1}\{z_1\}$ is a block,
it follows from the definition of a good solution that 
there is a natural one-to-one correspondence between the elements of $q^{-1}\{z_1\}$
and the blocks $\f A_i,$ $1\leq i \leq d.$ Furthermore, the action of any $\sigma\in G_h$ which maps $q^{-1}\{z_1\}$
to itself obviously respects this correspondence.  \qed

\section{Orbifolds on Riemann surfaces} 
In this section we recall main definitions and results related to Riemann surface orbifolds  (see  \cite{mil}, Appendix E). 

A pair $\f O=(R,\nu)$ consisting of a Riemann surface $R$ and a ramification function $\nu:R\rightarrow \mathbb N$ which takes the value $\nu(z)=1$ except at isolated points is called an orbifold. The Euler characteristic of an orbifold $\f O=(R,\nu)$ is defined by the formula  
$$ \chi(\f O)=\chi(R)+\sum_{z\in R}\left(\frac{1}{\nu(z)}-1\right),$$ where $\chi(R)$ is the 
Euler characteristic of $R.$

If $R_1$, $R_2$ are Riemann surfaces provided with ramification functions $\nu_1,$ $\nu_2$, and 
$f:\, R_1\rightarrow R_2$ is a holomorphic branched covering map, then $f$
is called a {\it covering map} $f:\,  \f O_1\rightarrow \f O_2$
{\it between orbifolds}
$\f O_1=(R_1,\nu_1)$ and $\f O_2=(R_2,\nu_2)$
if for any $z\in R_1$ the equality 
\be \l{us} \nu_{2}(f(z))=\nu_{1}(z)\deg_zf\ee holds, where $\deg_zf$ is the local degree of $f$ at the point $z$.
If for any $z\in R_1$ instead of equality \eqref{us} 
a weaker condition 
\be \l{uuss} \nu_{2}(f(z))\mid \nu_{1}(z)\deg_zf\ee
holds,  then $f$
is called a {\it holomorphic map} $f:\,  \f O_1\rightarrow \f O_2$
{\it between orbifolds}
$\f O_1$ and $\f O_2.$

A universal cover of an orbifold ${\f O}$
is a covering map between orbifolds \linebreak $\theta_{\f O}:\,
\tt {\f O}\rightarrow \f O$ such that $\tt R$ is simply connected and $\tt \nu(z)\equiv 1.$ 
If $\theta_{\f O}$ is such a map, then 
there exists a group $\Gamma_{\f O}$ of conformal automorphisms of $\tt R$ such that the equality 
$\theta_{\f O}(z_1)=\theta_{\f O}(z_2)$ holds for $z_1,z_2\in \tt R$ if and only if $z_1=\sigma(z_2)$ for some $\sigma\in \Gamma_{\f O}.$ A universal cover exists and 
is unique up to a conformal isomorphism of $\tt R,$
unless $\f O$ is the Riemann sphere with one ramified point, or $\f O$ is the Riemann sphere with two ramified points for which $\nu(z_1)\neq \nu(z_2).$ Unless stated otherwise, below we will assume that  considered orbifolds have a universal cover. Abusing  notation we will use the symbol $\tt {\f O}$ both for the
orbifold and for the  Riemann surface  $\tt R$.

If $f:\,  \f O_1\rightarrow \f O_2$ is a covering map 
between orbifolds, then for any choice of $\theta_{\f O_1}$ and $\theta_{\f O_2}$ there exists a conformal isomorphism $\tau:\, 
\tt{\f O_1}\rightarrow \tt{\f O_2}$ such that the diagram  
\be \l{edcv} 
\begin{CD}
\tt {\f O_1} @>\tau>> \tt {\f O_2}\\
@VV\theta_{\f O_1}V @VV\theta_{\f O_2}V\\ 
\f O_1 @>f >> \f O_2\ 
\end{CD}
\ee
is commutative. 
More generally, the following proposition holds.

\bp \l{poiu} Let $f:\,  \f O_1\rightarrow \f O_2$ be a holomorphic map between orbifolds. Then for any choice of $\theta_{\f O_1}$ and $\theta_{\f O_2}$ there exist 
a holomorphic map $F:\, \tt {\f O_1} \rightarrow \tt {\f O_2}$ and 
a homomorphism $\phi:\, \Gamma_{\f O_1}\rightarrow \Gamma_{\f O_2}$ such that diagram 
\be \l{dia2}
\begin{CD}
\tt {\f O_1} @>F>> \tt {\f O_2}\\
@VV\theta_{\f O_1}V @VV\theta_{\f O_2}V\\ 
\f O_1 @>f >> \f O_2\ 
\end{CD}
\ee
is commutative and 
for any $\sigma\in \Gamma_{\f O_1}$ the equality
\be \l{homm}  F\circ\sigma=\phi(\sigma)\circ F \ee holds.
The map $F$ is defined by $\theta_{\f O_1}$, $\theta_{\f O_2}$, and $f$  
uniquely up to the transformation 
$F\rightarrow g\circ F,$ where $g\in \Gamma_{\f O_2}$. 
In the other direction, for any holomorphic map $F:\, \tt {\f O_1} \rightarrow \tt {\f O_2}$ which satisfies \eqref{homm} for some homomorphism $\phi:\, \Gamma_{\f O_1}\rightarrow \Gamma_{\f O_2}$
there exists a uniquely defined  holomorphic map between orbifolds $f:\,  \f O_1\rightarrow \f O_2$ such that diagram \eqref{dia2} is commutative.
The holomorphic map $F$ is an isomorphism if and only if $f$ is a covering map between orbifolds.

\ep
\pr Assume that \eqref{uuss} holds. 
Let $F$ be the complete analytic continuation of
the germ 
$\theta_{\f O_2}^{-1}\circ f\circ \theta_{\f O_1}$, where 
$\theta_{\f O_2}^{-1}$ is a germ of a branch of the function inverse to $\theta_{\f O_2}$. 
Clearly, the local multiplicity of the map $f\circ \theta_{\f O_1}$ 
at a point $z\in \tt{\f O_1}$ equals $$\deg_z\theta_{\f O_1}\,\deg_{ \theta_{\f O_1}(z)}f=\nu(\theta_{\f O_1}(z))\deg_{ \theta_{\f O_1}(z)}f.
$$ On the other hand, the order of the permutation of branches of the function inverse to $\theta_{\f O_2}$ induced by the analytic continuation along a small loop around the point $(f\circ \theta_{\f O_1})(z)$ is equal to $\nu_2((f\circ \theta_{\f O_1})(z)).$ Therefore, condition  \eqref{uuss} 
implies that the function $F$ has no ramification points.
Since $\tt{\f O_1}$ is simply connected, we conclude that $F$ is single valued. 

Furthermore, if $\widehat F$ is another function which makes diagram \eqref{dia2} commutative, then $\theta_{\f O_2} \circ F=\theta_{\f O_2}\circ \widehat F$. 
Thus, for any $z\in \t R_1$ we have $\widehat F(z)=g_z\circ F(z)$ for some $g_z\in  \Gamma_{\f O_2}$. Since $\Gamma_{\f O_2}$ is countable this implies easily that 
$\widehat F \equiv g\circ F$ for some $g\in \Gamma_{\f O_2}$. Finally, since \eqref{dia2} implies that 
for any $\sigma\in \Gamma_{\f O_1}$ the equality $$\theta_{\f O_2} \circ F=\theta_{\f O_2}\circ (F\circ \sigma)$$ holds,  
we have: $$F\circ \sigma=\phi(\sigma)\circ F $$ for some $\phi(\sigma)\in \Gamma_{\f O_2}$, and it is easy to see that the correspondence $\sigma\rightarrow \phi(\sigma)$ 
is a homomorphism.

In the other direction, if \eqref{homm} holds, then $F$ maps any orbit of $\Gamma_{\f O_1}$  to an orbit of $\Gamma_{\f O_2}$, implying that the function $f=\theta_{\f O_2}\circ F\circ \theta_{\f O_1}^{-1}$ is well defined and holomorphic.
Further, \e{dia2} implies that 
$$\deg_z(f\circ \theta_{\f O_1})=\deg_z(\theta_{\f O_2}\circ F).$$ 
Therefore, since
$$\deg_z(f\circ \theta_{\f O_1})=\deg_z\theta_{\f O_1}\,\deg_{\theta_{\f O_1}(z)}f=\nu_1(\theta_{\f O_1}(z))\deg_{ \theta_{\f O_1}(z)}f$$ and 
\begin{multline*}  
\deg_z(\theta_{\f O_2}\circ F)=
\deg_zF\,\deg_{F(z)} \theta_{\f O_2}= \\
=\deg_zF\,\nu_2((\theta_{\f O_2}\circ F)(z))=\deg_zF\,\nu_2((f\circ \theta_{\f O_1})(z)),
\end{multline*}
we conclude that 
\be \l{wsx} \nu_1(\theta_{\f O_1}(z))\deg_{ \theta_{\f O_1}(z)}f=\deg_zF\,\nu_2((f\circ \theta_{\f O_1})(z)),\ee
implying  
\eqref{uuss}. 
Moreover, it follows from \eqref{wsx} that $F$ is locally and therefore globally invertible if and only if \eqref{us} 
holds.
\qed
\vskip 0.2cm

If $f:\,  \f O_1\rightarrow \f O_2$ is a covering map between orbifolds with compact support, then  the Riemann-Hurwitz 
formula implies that 
\be \l{rhor} \chi(\f O_1)=d \chi(\f O_2), \ee
where $d=\deg f$. 
More generally, the following statement is true. 

\bp \l{p22} Let $f:\, \f O_1\rightarrow \f O_2$ be a holomorphic map between orbifolds with compact support.
Then 
\be \l{iioopp} \chi(\f O_1)\leq \chi(\f O_2)\,\deg f, \ee and the equality 
holds if and only if $f:\, \f O_1\rightarrow \f O_2$ is a covering map between orbifolds.

\ep 
\pr Denote by $S_1$ (resp. $S_2$) the set of ramified points of $\f O_1$ (resp.  $\f O_2$) and by $\f C(f)$ the set of critical values of
$f.$ Set $$S=f(S_1)\cup \f C(f),  \ \ \ \widehat R_2=R_2\setminus S, \ \  \  \widehat R_1=f^{-1}\{\widehat R_2\}.$$ 
Observe that $S_2\subseteq S$ since \eqref{uuss} implies that whenever $\nu_2(f(z))$ is greater than one at least one of the numbers $\deg_zf$ and $\nu_1(z)$ also is greater than one.   
Since  $f:\, \widehat R_1 \rightarrow \widehat R_2$ is a covering map between surfaces, we have: $$\chi(\widehat R_1)=d\chi(\widehat R_2),$$ where $d=\deg f.$ 
Furthermore, it follows from \eqref{uuss} that  
$$\frac{1}{\nu_{1}(z)}\leq \frac{\deg_zf}{\nu_{2}(f(z))}$$ implying that 
\be \l{ppqq} \sum_{\substack{x\in R_1\\ f(x)=f(z)}}\frac{1}{\nu_{1}(x)}\leq \frac{d}{\nu_{2}(f(z))},\ee where the equality holds if and only if \e{us} 
holds for any $x\in f^{-1}\{z\}.$

Since removing a point from a surface reduces the Euler characteristic by one, we have:
\begin{multline*}  
\chi(\f O_1)=  \chi(R_1)+\sum_{\substack{x\in R_1}}\left(\frac{1}{\nu_1(x)}-1\right)=\chi(R_1)+\sum_{\substack{x\in R_1\\ f(x)\in S}}\left(\frac{1}{\nu_1(x)}-1\right)=\\ 
=\chi(\widehat R_1)+\sum_{\substack{x\in R_1\\ f(x)\in S}}\frac{1}{\nu_1(x)}=d\chi(\widehat R_2)+\sum_{\substack{x\in R_1\\ f(x)\in S}}\frac{1}{\nu_1(x)}.
\end{multline*}  
It follows now from \eqref{ppqq} that 
$$\chi(\f O_1)\leq d\chi(\widehat R_2)+ \sum_{z \in S}\frac{d}{\nu_2(z)}=d\chi( \f O_2),$$
where the equality holds if and only if \e{us} 
holds for any $z\in f^{-1}\{S\}.$
Since the definition of $S$ implies that \e{us} is satisfied for 
$z\notin f^{-1}\{S\}$, we conclude that the equality in \e{iioopp} holds if and only if 
$f\,:\, \f O_1\rightarrow \f O_2$ is a covering map between orbifolds. \qed

\bc \l{mc} 
Let $f:\,\f O_1\rightarrow \f O_2$ be a holomorphic map  between orbifolds  with compact supports. Assume that $\deg f>1$ and $\chi(\f O_1)=\chi(\f O_2)=l.$ Then $l \geq 0$. Furthermore, $l=0$ if and only if  $f$ is a covering map. \qed
\ec

\section{
Minimal maps and decompositions}
In this section we introduce the concept of a minimal holomorphic map between orbifolds, and 
establish some properties of such maps related to functional decompositions.

Let $R_1$, $R_2$ be Riemann surfaces, and 
$f:\, R_1\rightarrow R_2$ a holomorphic branched covering map. Assume that $R_2$ is provided with a ramification function $\nu_2$. In order to define a ramification function $\nu_1$ on $R_1$ so that $f$ would be a holomorphic map between orbifolds $\f O_1=(R_1,\nu_1)$ and $\f O_2=(R_2,\nu_2)$ 
we must satisfy condition \eqref{uuss}, and it is easy to see that
for any  $z\in R_1$ a minimal possible value for $\nu_1(z)$ is defined by 
the equality 
\be \l{rys} \nu_{2}(f(z))=\nu_{1}(z)\GCD(\deg_zf, \nu_{2}(f(z)).\ee 
In case if \eqref{rys} is satisfied for  any $z\in R_1$ we 
say that $f$ is a {\it minimal holomorphic  map 
between orbifolds} 
$\f O_1=(R_1,\nu_1)$ and $\f O_2=(R_2,\nu_2)$.

It follows from the definition that for any orbifold $\f O=(R,\nu)$ and holomorphic branched covering map $f:\, R^{\prime} \rightarrow R$ there exists a {\it unique} orbifold structure $\nu^{\prime}$ on $R^{\prime}$ such that 
$f$ is a minimal holomorphic map between corresponding orbifolds. We will denote the corresponding orbifold by $f^*\f O.$

\bl \l{wert} 
Any covering map between orbifolds $f:\, \f O_1\rightarrow \f O_2$  is a  minimal holomorphic map. In particular, the equality $\f O_1= f^*\f O_2$ holds. A  minimal holomorphic map $f:\, \f O_1\rightarrow \f O_2$  is 
a covering map if and only if $\deg_zf\mid \nu_2(f(z))$ for any $z\in R_1.$
\el
\pr Follows from the corresponding  definitions. \qed 

\vskip 0.2cm
Minimal holomorphic maps possess the following fundamental property with respect to compositions.

\bt \l{41} Let $f:\, R^{\prime\prime} \rightarrow R^{\prime}$ and $g:\, R^{\prime} \rightarrow R$ be holomorphic branched covering maps, and  $\f O=(R,\nu)$ an orbifold. 
Then 

$$(g\circ f)^*\f O= f^*(g^*\f O).$$
\et 
\pr 
Let $f^*(g^*\f O)=(R^{\prime\prime},\nu_1)$ and $g^*\f O=(R^{\prime},\nu_2)$. 
Since $$f:\, f^*(g^*\f O)\rightarrow g^*\f O$$
and 
$$g:\, g^*\f O\rightarrow \f O$$
are minimal holomorphic maps,
for any $z\in  R^{\prime\prime}$ we have:
\be \l{11} \nu_{2}(f(z))=\nu_{1}(z)\GCD\Big(\deg_zf, \nu_{2}(f(z))\Big)\ee
and \be \l{12} \nu\Big((g\circ f)(z)\Big)=\nu_{2}(f(z))\GCD\bigg(\deg_{f(z)}g, \nu\Big((g\circ f)(z)\Big)\bigg).\ee 
In order to prove the theorem, we only must show that
\be \l{13}  \nu\Big((g\circ f)(z)\Big)=\nu_{1}(z)\GCD\bigg(\deg_{z}(g\circ f), \nu\Big((g\circ f)(z)\Big)\bigg).
\ee

Observe first that for any positive integers $a,b,c$ the equality
\be \l{00} \GCD(ab,c)=\GCD(a,c)\GCD\left(b,\frac{c}{\GCD(a,c)}\right)\ee
holds. Indeed, the last statement is equivalent to the statement  that for any non-negative integers $\alpha,\beta,\gamma$  the equality 
\be \l{001} 
\min\{\alpha+\beta, \gamma\}=\min\{\alpha, \gamma\}+
\min\{\beta, \gamma-\min\{\alpha, \gamma\}\}\ee
holds. If $\min\{\alpha, \gamma\}= \gamma,$ then clearly  $\gamma\leq \alpha+\beta$ and inequality \eqref{001} is true.
On the other hand, if  $\min\{\alpha, \gamma\}= \alpha,$ then \eqref{001} reduces to 
to the obvious equality 
$$
\min\{\alpha+\beta, \gamma\}=\alpha+
\min\{\beta, \gamma-\alpha\}.$$

Setting  in \eqref{00}
$$a=\deg_{f(z)}g, \ \ \ b=\deg_zf, \ \ \ c= \nu\Big((g\circ f)(z)\Big).$$ 
and using the formulas \be \l{000} \deg_{z}(g\circ f)=\deg_zf\,\deg_{f(z)}g\ee
and \eqref{12}, we have: 
\begin{multline} \l{kjh} 
\GCD\bigg(\deg_{z}(g\circ f), \nu\Big((g\circ f)(z)\Big)\bigg)= \\
=\GCD\bigg(\deg_{f(z)}g, \nu\Big((g\circ f)(z)\Big)\bigg)\GCD\Big(\deg_zf, \nu_{2}(f(z))\Big). 
\end{multline}
Since  equalities \eqref{11} and \eqref{12} 
imply the equality   
$$\nu\Big((g\circ f)(z)\Big)=\nu_{1}(z)\GCD\Big(\deg_zf, \nu_{2}(f(z))\Big)\GCD\bigg(\deg_{f(z)}g, \nu\Big((g\circ f)(z)\Big)\bigg)$$
equality \eqref{13} follows now from equality 
\eqref{kjh}. \qed

\bc \l{indu1} Let $f:\, \f O_1\rightarrow \f O^{\prime}$ and $g:\, \f O^{\prime}\rightarrow \f O_2$ be minimal holomorphic maps (resp. covering maps) between orbifolds.
Then  $g\circ f:\, \f O_1\rightarrow \f O_2$ is  a minimal holomorphic map (resp. a covering map). 
\ec 
\pr It follows from the conditions that 
$$\f O^{\prime}=g^* \f O_2, \ \ \  \f O_1=f^*(\f O^{\prime})=f^*(g^*\f O_2).$$
Therefore,  by Theorem \ref{41}, $$\f O_1=(g\circ f)^*\f O_2,$$ implying that $g\circ f:\, \f O_1\rightarrow \f O_2$  is  a minimal holomorphic map. 
Furthermore, the equalities  
\be \l{111} \nu^{\prime}(f(z))=\nu_{1}(z)\deg_zf,\ee
and \be \l{121} \nu_2\Big((g\circ f)(z)\Big)=\nu^{\prime}(f(z))\deg_{f(z)}g\ee 
imply the equality  
\be \l{123} \nu_2\Big((g\circ f)(z)\Big)=\nu_{1}(z)\deg_{z}(g\circ f). \ \ \ \Box\ee

\vskip 0.2cm

\bc \l{indu2}  Let $f:\, R_1 \rightarrow R^{\prime}$ and $g:\, R^{\prime} \rightarrow R_2$ be holomorphic branched covering maps, and  $\f O_1=(R_1,\nu_1)$ and  $\f O_2=(R_2,\nu_2)$
orbifolds. Assume that  $g\circ f:\, \f O_1\rightarrow \f O_2$ is  a minimal holomorphic map (resp. a co\-vering map). Then  $f:\, \f O_1\rightarrow g^*\f O_2 $ and $g:\, g^*\f O_2\rightarrow \f O_2$ are minimal holomorphic maps (resp. covering maps). 
\ec 
\pr Set $\f O^{\prime}=g^* \f O_2.$ By condition, $\f O_1=(g\circ f)^*\f O_2$. Therefore,  by Theorem \ref{41}, $$\f O_1=f^*(\f O^{\prime}),$$ implying that $f:\, \f O_1\rightarrow  \f O^{\prime}$ is a minimal holomorphic map. 

Further, if $g\circ f:\, \f O_1\rightarrow \f O_2$ is a covering map, then it follows from 
\eqref{123} that 
%$$\deg_{z}(g\circ f)\mid \nu_2\Big((g\circ f)(z)\Big).$$ Therefore, 
$$\deg_{f(z)}g\mid \nu_2\Big((g\circ f)(z)\Big),$$ implying that $g:\, \f O^{\prime}\rightarrow \f O_2$ is a covering map.
Now, equalities \eqref{123} and \eqref{121} imply equality \eqref{111}. \qed

\vskip 0.2cm

With each holomorphic function $f:\, R_1\rightarrow R_2$ between compact Riemann surfaces 
one can associate in a natural way two orbifolds $\f O_1^f=(R_1,\nu_1^f)$ and 
$\f O_2^f=(R_2,\nu_2^f)$, setting $\nu_2^f(z)$  
equal to the least common multiple of local degrees of $f$ at the points 
of the preimage $f^{-1}\{z\}$, and $$\nu_1^1(z)=\frac{\nu_2^f(f(z))}{\deg_zf}.$$ By construction, 
$f$ is a covering map between orbifolds $f:\, \f O_1^f\rightarrow \f O_2^f.$ Notice that by Lemma \ref{wert} the equality  
\be \l{zzxx} \f O_1^f=f^* \f O_2^f\ee
holds. 

\bl \l{have} Orbifolds $\f O_1^f$ and 
$\f O_2^f$ have a universal cover. 
\el
\pr
Equality \e{posl} implies that $f$ may not have only one critical value. Moreover, if $f$ has two critical values, then, by transitivity 
of $G_h,$ the corresponding permutations have the same order equal to $\deg f$. Therefore,
$\f O_2^f$ has a universal cover.

Show now that  $\f O_1^f$ also has a universal cover.
Let $\theta:\, \tt{\f O}\rightarrow \f O_2^f$ be a universal cover of $\f O_2^f$,
and  $\widehat \theta$ the complete analytic continuation of a germ $f^{-1}\circ \theta,$ 
where $f^{-1}$ is a germ of a branch of the algebraic function inverse
to $f.$ It follows from the equality
\be \l{uuuii}  \deg_{z}\theta =\nu_{2}^f(\theta(z)), \ \ z\in \tt{\f O}, \ee
and the definition of $\f O_2^f$
%\be \l{uuuiii} 
%\ \ \ \ \ \nu_{2}^f(f(z))=\nu_{1}^f(z)\deg_{z}f, \ \ z\in \C\P^1 \ee
that $\widehat \theta$ has no local branching. Therefore, since $\tt{\f O}$ is simply connected, 
$\widehat \theta$ is single valued.
Moreover,  $f\circ \widehat \theta=\theta$. It follows now from equality  \eqref{zzxx} and Corollary \ref{indu2}  that $\widehat\theta:\, \tt{\f O}\rightarrow \f O_1^f$ is a covering map between orbifolds.
Since $\tt{\f O}$ is non-ramified, this implies that 
$\widehat\theta$ is 
a universal cover of $\f O_1^f.$ \qed
\vskip 0.2cm

\bt \l{t1} Let $h,f,p,g,q$ be a good solution of the equation $h=f\circ p=g\circ q.$ 
Then the commutative diagram 
\be 
\begin{CD}
\f O_1^q @>p>> \f O_1^f\\
@VV q V @VV f V\\ 
\f O_2^q @>g >> \f O_2^f\ 
\end{CD}
\ee
consists of minimal holomorphic  maps between orbifolds

\et 

\pr 
Denote by $C$ the Riemann surface on which $g$ is defined.
Let $z \in C$ be a point, $\rho\subset C$ a small free loop around $z$, and 
$z_1\in \rho $ a point such that $g(z_1)=z_0$ is a regular value of $h$. Then a permutation of points of $h^{-1}\{z_0\}$ corresponding to the 
analytic continuation of $h^{-1}$ along the curve $g(\rho)\subset \C\P^1$ induces a permutation $\sigma_1$
of points of $q^{-1}\{z_1\}$ as well as a permutation $\sigma_2$ of points 
of $f^{-1}\{z_0\}$. Furthermore, by Lemma \ref{good1} 
the permutations $\sigma_1$ and $\sigma_2$ have the same cyclic structure. 
In particular, since the order of $\sigma_1$ is equal to  $\nu_2^{q}(z)$ by construction of $\f O_2^q,$  the order of $\sigma_2$ also is equal to $\nu_2^{q}(z)$. 

On the other hand, $$\sigma_2=\sigma_3^{\deg_zg},$$ where 
$\sigma_3$ is a permutation
of points $f^{-1}\{z_0\}$ induced by the analytic continuation 
of $h^{-1}$ along a small free loop $\tt \rho$ around $g(z).$ Since the order of $\sigma_3$ is equal to $\nu_2^{f}(g(z))$, this yields that 
the order of $\sigma_2$ is equal to  
$$\frac{\nu_2^{f}(g(z))}{\GCD(\deg_zg, \nu_2^{f}(g(z))},$$ implying the equality 
\be \l{orr} \nu_2^{f}(g(z))=\nu_2^{q}(z)\GCD(\deg_zg, \nu_2^{f}(g(z)).\ee
Thus, $g:\, {\f O_2^q}\rightarrow {\f O_2^f}$ is a minimal holomorphic  map between orbifolds.

Further, since $g:\, {\f O_2^q}\rightarrow {\f O_2^f}$ is a minimal holomorphic  map, 
it follows from Corollary \ref{indu1} that $g\circ q :\, {\f O_1^q}\rightarrow {\f O_2^f}$ also is a minimal holomorphic  map.
Finally, since $f\circ p=g\circ q$, it follows now from   Corollary \ref{indu2} taking into account equality \eqref{zzxx} that $p:\, {\f O_1^q}\rightarrow {\f O_1^f}$ also
is a minimal holomorphic  map. \qed

\section{Minimal self-maps and equivariant functions}
Let $f$ be a rational function of degree at least two such that  $f\,: \f O\rightarrow \f O$ is a minimal holomorphic self-map between orbifolds for some $\f O$ defined on $\C\P^1.$
Then by Corollary \ref{mc} the inequality $\chi(\f O)\geq 0$ holds and the equality attains if and only if $f$ is a covering map. In  turn, the condition that $f\,: \f O\rightarrow \f O$ is a covering map for some 
$\f O$ with $\chi(\f O)= 0$ is equivalent to the condition that $f$ is a Latt\`es function (see \cite{mil2}, Theorem 4.1). Thus,  minimal 
holomorphic self-maps $f\,: \f O\rightarrow \f O$ for $\f O$ with $\chi(\f O)= 0$
are exactly Latt\`es functions whose properties are well-established 
(see \cite{mil2} for a   comprehensive survey of these properties). In contrast, rational functions which are minimal holomorphic self-maps $f:\, \f O\rightarrow \f O$ for $\f O$ with $\chi(\f O)>0$ seem to be a completely  new object. In this section we  provide a characterization  of such functions, and outline an approach to their classification.

%It is well known  that 
If $\chi(\f O)>0$ and $\f O$
is neither non-ramified sphere nor one of two orbifolds without the universal cover, then the collection of ramification indices of $\f O$ is either $(n,n)$, or $(2,2,n)$ for some $n\geq 2,$ or one of the following triples $(2,3,3)$, $(2,3,4),$ $(2,3,5).$ 
Further, for all these collections the universal cover $\tt{\f O}$ of $\f O$ is $\C\P^1$, and the corresponding groups $\Gamma_{\f O}$ 
are finite subgroups of the automorphism group of $\C\P^1$, 
namely, the cyclic, dihedral, tetrahedral, octahedral, and icosahedral. Accordingly, the functions $\theta_{\f O}$ are rational functions of degree $n,$ $2n$, $12,$ $24$, and $60$ which can be characterized as  regular covers of the sphere by the sphere
(see e.g. \cite{k}).

\bt \l{las}  Let $A$ and $F$ be rational functions of degree at least two and  $\f O=(\C\P^1, \nu)$ an orbifold with $\chi(\f O)>0$ such that 
$A:\, \f O \rightarrow \f O$ 
is a holomorphic map between orbifolds  and  the diagram 
\be \l{dia3}
\begin{CD}
\tt {\f O} @>F>> \tt {\f O}\\
@VV\theta_{\f O}V @VV\theta_{\f O}V\\ 
\f O @>A >> \f O\ 
\end{CD}
\ee
is commutative.
Then the following conditions are equivalent.

\begin{enumerate}
\item The holomorphic map $A$ is a minimal holomorphic  map. 
\item  The homomorphism $\phi:\, \Gamma_{\f O}\rightarrow \Gamma_{\f O}$ defined by the equality 
\be \l{homo} F\circ\sigma=\phi(\sigma)\circ F, \ \ \ \sigma\in \Gamma_{\f O},\ee is an automorphism of $\Gamma_{\f O}$.
\item The triple $F,$ $A,$ $\theta_{\f O}$ is a good solution of the equation 
\be \l{qws} 
A\circ \theta_{\f O}=\theta_{\f O}\circ F.\ee 
\end{enumerate}

\et

\pr By Lemma \ref{good}, in order to prove $2\Leftrightarrow 3$ it is enough to show that $\phi$ is an automorphism 
if and only if the functions 
$F$ and $\theta_{\f O}$ have no non-trivial common compositional right factor. 
Since $\theta_{\f O}$ is a regular cover, any compositional right factor of  $\theta_{\f O}$ of degree greater than one has the form $\theta_{\f O^{\prime}}$, 
where  $\Gamma_{\f O^{\prime}}\neq \{e\}$ is a subgroup  of $\Gamma_{\f O}$.
Therefore, the functions $F$ and $\theta_{\f O}$ have a non-trivial common compositional right factor if and only if there exists $\Gamma_{\f O^{\prime}}\subseteq \Gamma_{\f O}$ such that 
$F(z_1)=F(z_2)$ whenever $z_2=\sigma(z_1)$ for some $\sigma \in \Gamma_{\f O^{\prime}}$. On the other hand, by \eqref{homo}, such a subgroup exists if and only of $\phi$  is not a monomorphism. This proves the equivalence  $2\Leftrightarrow 3$.

By Lemma \ref{good} and Corollary \ref{good2}, in order to prove $1\Leftrightarrow 3$, 
it is enough to show that 
$A:\, {\f O}\rightarrow {\f O}$ is a minimal holomorphic  map between orbifolds if and only if $\deg_z\theta_{\f O}$ and 
$\deg_zF$ are coprime for any $z\in \C\P^1.$ Since 
$$\deg_z \theta_{\f O} \deg_{\theta_{\f O}(z)}A=\deg_zF\deg_{F(z)}\theta_{\f O},$$
the last
condition is equivalent to the equality
\be \l{asz} \deg_{F(z)}\theta_{\f O}=
\deg_z \theta_{\f O} \GCD(\deg_{\theta_{\f O}(z)}A,\deg_{F(z)}\theta_{\f O}).\ee
On the other hand, since $\deg_z{\theta_{\f O}}=\nu(\theta_{\f O}(z))$ and 
$$\deg_{F(z)}\theta_{\f O}=\nu(\theta_{\f O}(F(z)))=\nu(A(\theta_{\f O}(z))),$$ equality \eqref{asz} is equivalent to the equality 
$$\nu(A(\theta_{\f O}(z)))=
 \nu(\theta_{\f O}(z))\GCD\Big(\deg_{\theta_{\f O}(z)}A,\nu(A(\theta_{\f O}(z)))\Big),$$ 
which in  turn is equivalent
to the 
requirement that $A:\, {\f O}\rightarrow {\f O}$ 
is a minimal holomorphic  map. \qed

\vskip 0.2cm
\bc \l{epta} Let $\f O=(\C\P^1, \nu)$ be an orbifold with $\chi(\f O)>0,$ and $A:\, \f O \rightarrow \f O$ 
a minimal holomorphic map between orbifolds.  Then   for any decomposition $A=U\circ V$ the ramification collection of  $U^*\f O$ coincides with the one of $\f O$.
\ec
\pr Denote $U^*\f O$  by $\f O^{\prime}.$ It follows from Corollary \ref{indu2} and  Proposition \ref{poiu} that there 
exist holomorphic maps $F_U$ and $F_V$ which make the diagram    
\begin{equation*} 
\begin{CD} 
\tt{\f O} @>F_V>> \tt{\f O^{\prime}} @>F_U>>\tt{\f O} \\
@VV\theta_{\f O}  V @VV \theta_{\f O^{\prime}} V @VV \theta_{\f O} V\\ 
\f O @> V >> \f O^{\prime} @>U>>\f O \\
\end{CD}
\end{equation*}
commutative. 
Since $$\chi(\f O)\leq \chi(\f O^{\prime})\deg V$$ by Proposition \ref{p22}, the inequality $\chi(\f O^{\prime})>0$ holds, implying that  
$\tt{\f O^{\prime}}=\C\P^1$. 
Furthermore, there  exist homomorphisms 
$$\phi_V:   \Gamma_{\f O}\rightarrow \Gamma_{\f O^{\prime}}, \ \ \  \phi_U:   \Gamma_{\f O^{\prime}}\rightarrow \Gamma_{\f O}$$
such that  
$$F_V\circ\sigma=\phi_V(\sigma)\circ F_V, \ \ \sigma\in \Gamma_{\f O}, \ \ \ \ 
F_U\circ\sigma=\phi_U(\sigma)\circ F_U, \ \ \ \sigma\in \Gamma_{\f O^{\prime}}.$$
Since the function $F_U\circ F_V$ makes diagram \eqref{dia3} commutative, 
 Theorem \ref{las} implies that the composition of homomorphisms $$\phi_{U}\circ \phi_{V}: \Gamma_{\f O}\rightarrow \Gamma_{\f O}$$ is an 
automorphism. Therefore, $\Gamma_{\f O^{\prime}}\cong\Gamma_{\f O}$, and the ramification collection of  $\f O^{\prime}$ coincides with the one of $\f O$. \qed

\vskip 0.2cm
 Theorem \ref{las} reduces the study of minimal holomorphic maps $f:\, \f O\rightarrow \f O$ with $\chi(\f O)>0$ to the study of rational functions 
such that \eqref{homo} holds for some finite subgroup $\Gamma_{\f O}$ of $\C\P^1$ and an automorphism $\phi:\, \Gamma_{\f O}\rightarrow \Gamma_{\f O}$.
For example, consider an orbifold $\f O$ with ramification $$\nu(0)=n, \ \ \ \nu(\infty)=n.$$ 
In this case,
$\Gamma_{\f O}=\Z/n\Z$ is generated by the transformation $$E_n:\, z\rightarrow e^{2\pi i/n} z,$$ and $\theta_{\f O}=z^n.$
 Since a homomorphism $\phi: \Gamma_{\f O}\rightarrow  \Gamma_{\f O}$ is an automorphism if and only if 
\be \l{transa2} \phi(E_n): z\rightarrow e^{2\pi ir/n} z,\ee for some $r,$ $1\leq r\leq n-1,$ with $\GCD(r,n)=1$, a rational function $F$ satisfies \eqref{homo}  for some automorphism $\phi$ 
of $\Gamma_{\f O}$
if and only if $F/z^r$ is $\Gamma_{\f O}$-invariant
for some  $r$ as above, that is if and only if 
$$F=z^rR(z^n)$$ for some rational function $R$. 
Accordingly, since the function $A=z^rR^n(z)$ makes corresponding diagram \eqref{dia3} commutative,   minimal holomorphic maps $f:\, \f O\rightarrow \f O$ have the form 
$z^rR^{n}(z),$ and 
the correspondence between  minimal holomorphic maps and functions satisfying 
\eqref{homo} is given by 
the commutative  diagram 
\be \l{zaqw} 
\begin{CD}
\tt {\f O} @>z^rR(z^n)>> \tt {\f O}\\
@VVz^nV @VV z^nV\\ 
\f O @>z^rR^n(z) >> \f O.\ 
\end{CD}
\ee

More generally, it was shown in \cite{dm} that for finite subgroups $\Gamma$ of $Aut(\C\P^1)$ 
the problem of describing of $\Gamma$-equivariant functions, that is of functions $F$ satisfying \eqref{homo} for the {\it identity} automorphism $\phi$, may be reduced to the classical problem of describing of homogeneous $\Gamma$-invariant polynomials solved by Klein in \cite{k}. 
For the group $\Gamma=S_4$, say, this solves the problem. Indeed, since $F$ corresponding to $A$ in \eqref{dia3} is defined up to the change 
$F\rightarrow g\circ F$, where $g\in \Gamma_{\f O}$  (see Proposition \ref{poiu}), the automorphism $\phi$ is defined up to the change
$\phi\rightarrow g\circ \phi \circ g^{-1}$. Therefore, since all automorphisms of $S_4$  are inner, without loss of generality we may assume that $\phi(\sigma)=\sigma$. 
For $\Gamma$ equal to $A_4$ or $A_5$ 
the group  $Out(\Gamma)$ consists of two elements, implying that in these cases
in addition to the  identity 
automorphism we must consider one additional automorphism,  and  one can try to extend
the approach of \cite{dm} so that to cover this case too. Notice also that although a function $F$ satisfying \eqref{homo} is not necessary  $\Gamma_{\f O}$-equivariant, it follows from the finiteness of  $\Gamma_{\f O}$ that  some {\it iterate} of $F$ is  $\Gamma_{\f O}$-equivariant.

\section{Proof of Theorem \ref{main} and examples}
We start from proving Theorem \ref{main} in the case where $\C(B,X)=\C(z)$. We will call solutions of \eqref{i1} satisfying this condition {\it primitive}.
By Lemma \ref{good}, a solution  $A,X,B$ of \eqref{i1} is primitive if and only 
if the corresponding solution
of \e{m}  given by   
$$f=q=X, \ \ p=B, \ \ g=A,$$ is good.   
Clearly, any solution $A,X,B$ of \eqref{i1} with $\deg X=1$ is primitive. The corresponding functions $A$ and $B$ are conjugated and therefore equivalent.

\bt \l{oip} Let  $A,X,B$ be a primitive solution of \e{i1} with $\deg X>1$. Then $\chi(\f O_1^X)\geq 0$, $\chi(\f O_2^X)\geq 0$, and 
the commutative diagram 
\be 
\begin{CD} \l{gooopa}
\f O_1^X @>B>> \f O_1^X\\
@VV X V @VV X V\\ 
\f O_2^X @>A >> \f O_2^X\ 
\end{CD}
\ee
consists of minimal holomorphic  maps between orbifolds.  Furthermore, either $\chi(\f O_1^X)= \chi(\f O_2^X)=0$
and $A$, $B$ are Latt\`es functions, or $0<\chi(\f O_2^X)< \chi(\f O_1^X)$. In the last case, possible  
collections of ramification indices of $\f O_1^X$ and $\f O_2^X$ are following: $(n,n)$ or $(2,2,n)$ for some $n\geq 2,$ or one of the triples $(2,3,3)$, $(2,3,4),$ $(2,3,5).$ 
In addition, $\f O_1^X$ may be a non-ramified sphere. 
\et

\pr 
Since  $A:\, \f O_2^X\rightarrow \f O_2^X$ and $B:\, \f O_1^X\rightarrow \f O_1^X$ are minimal holomorphic  maps between orbifolds  by Theorem \ref{t1}, it follows from Corollary \ref{mc}  that $\chi(\f O_1^X)\geq 0$ and $\chi(\f O_2^X)\geq 0$. Furthermore, since $X:\, \f O_1^X\rightarrow \f O_2^X$ is a covering map, formula \eqref{rhor} implies that 
either  $\chi(\f O_1^X)= \chi(\f O_2^X)=0$, or $0<\chi(\f O_2^X)< \chi(\f O_1^X)$. In the first case,  $A$ and $B$ are Latt\`es functions.
On the other hand, if $\chi(\f O_2^X)>0$, then the collections of ramification indices of $\f O_1^X$ and $\f O_2^X$ have the required form since $\f O_1^X$ and $\f O_2^X$ have a universal cover by Lemma \ref{have},
and the condition $\deg X>1$ implies that $\f O_2^X$ is distinct from the non-ramified sphere. 
\qed

\vskip 0.2cm

The proof of Theorem \ref{main} in the general case reduces to the primitive case as follows.
Let $A,X,B$ be a non-primitive solution of \e{i1} and $U_1$ any rational function such that \be \l{tog0} \C(X,B)=\C(U_1)\ee so that 
\be \l{tog} X=X_1\circ U_1, \ \ \ B=V_1\circ U_1\ee for some $X_1,V_1\in \C(z)$ such that 
$\C(X_1,V_1)=\C(z).$ Clearly, \eqref{i1} and \eqref{tog} imply that 
$$A\circ X_1=X_1\circ (U_1\circ V_1).$$ Thus, $A,X_1,U_1\circ V_1$ is another solution of 
\eqref{i1}. Continuing in this way, 
define $X_{i+1},V_{i+1},$ and $U_{i+1}$, $i\geq 1$,     as rational functions satisfying  the equalities   
\be \l{tog1} \C(X_i,U_i\circ V_i)=\C(U_{i+1}),\ee  
\be \l{tog2} X_i=X_{i+1}\circ U_{i+1},\ee
\be \l{tog3}  U_i\circ V_i=V_{i+1}\circ U_{i+1},\ee
and set $$B_i=U_i\circ V_i.$$
Clearly, 
$$A\circ X_i=X_i\circ B_i.$$ Furthermore,  since 
$\deg X_{i+1}< \deg X_{i}$ whenever $\deg U_{i+1}>1,$ the equality $\C(X_l,B_l)=\C(z)$ holds for some $l\geq 1$, and hence   
$A,$ $X_l,$ $B_l$ is a primitive solution of \eqref{i1}.  

By construction, if $\deg X_l=1,$ then $B\sim A$. Otherwise, by Theorem \ref{oip}, 
the commutative diagram 
$$ 
\begin{CD} 
\f O_1^{X_l} @>B_l>> \f O_1^{X_l}\\
@VV X_l V @VV X_l V\\ 
\f O_2^{X_l} @>A >> \f O_2^{X_l}\ 
\end{CD}
$$ consists of minimal holomorphic  maps between orbifolds. 
Set 
$$\f O=U_l^* \f O_1^{X_l}.$$ Since $B_l=U_l\circ V_l,$ it follows from Corollary \ref{indu2} and 
Corollary \ref{epta} 
that $$V_l:\f O_1^{X_l}\rightarrow \f O, \ \ \  U_l:\f O\rightarrow \f O_1^{X_l}$$ are
minimal holomorphic  maps between orbifolds, and  $\f O$ has the same ramification collection  as $\f O_1^{X_l}$ (and hence the same Euler characteristic).
 It follows now from Corollary \ref{indu1}, taking into account the equality $$B_{l-1}
=U_{l-1}\circ V_{l-1}=V_l\circ U_l,$$ that
$$X_l\circ U_l: \f O\rightarrow \f O_2^{X_l}, \ \ \ B_{l-1}: \f O\rightarrow \f O$$ 
are minimal holomorphic maps between orbifolds.  
Since 
\be \l{xori+} X=X_l\circ U_l\circ U_{l-1}\circ \dots \circ U_1,\ee continuing in this way we see that the conclusion of Theorem \ref{main} holds for $\f O_2= \f O_2^{X_l}$ 
and some orbifold $\f O_1$ whose ramification collection coincides with the one of  $\f O_1^{X_l}.$ \qed

\vskip 0.2cm

The simplest examples of solutions of \eqref{i1}  are obtained
from diagram \eqref{dia3}, where $F$, $A$, and $\f O$ satisfy the conditions of Theorem \ref{las}. 
For example,  diagram \eqref{zaqw} provides a family of such examples. 
Moreover, since for any finite subgroup $\Gamma$ of $Aut(\C\P^1)$ there exist $\Gamma$-equivariant rational functions $F$ (see \cite{dm}), 
similar examples can be given for  any finite subgroup $\Gamma$ of $Aut(\C\P^1)$.

More generally, assume that the automorphism $\phi$ in \eqref{homo} satisfies the equality
$$\phi(\Gamma_{\f O^{\prime}})= \Gamma_{\f O^{\prime}}$$ for some subgroup $\Gamma_{\f O^{\prime}}$ of $\Gamma_{\f O}$ (notice that if $F$ is  $\Gamma_{\f O}$-equivariant this is true for any subgroup
$\Gamma_{\f O^{\prime}}$ of $\Gamma_{\f O}$).
Then by Proposition \ref{poiu} there exists a rational function 
$B$ such that \be \l{tre} B\circ \theta_{\f O^{\prime}}=\theta_{\f O^{\prime}}\circ F.\ee
On the other hand, it follows from $\Gamma_{\f O^{\prime}}\subset\Gamma_{\f O}$  that \be \l{iuy} \theta_{\f O}=X\circ \theta_{\f O^{\prime}}\ee for some rational function $X$.
Thus,  
the  diagram 
\be 
\begin{CD} \l{gooopa2}
\tt{\f O^{\prime}} @>F>> \tt{\f O^{\prime}}\\
@VV \theta_{\f O^{\prime}} V @VV \theta_{\f O^{\prime}} V\\ 
\f O^{\prime} @>B>> \f O^{\prime}\\
@VV X V @VV X V\\ 
\f O @>A >> \f O\ 
\end{CD}
\ee
is commutative, implying that $A,X,$ and $B$ satisfy \eqref{i1}.

In order to illustrate the above construction, consider an orbifold $\f O$  with $$\nu(1)=\nu(-1)=2, \ \ \ \nu(\infty)=n.$$ In this case $\Gamma_{\f O}=D_{2n}$
is generated by 
the transformations $$\alpha:\, z\rightarrow  e^{2\pi i/n} z, \ \ \ \beta:\, z\rightarrow \frac{1}{z},$$  and
\be \l{xoi} \theta_{\f O}=\frac{1}{2}\left(z^n+\frac{1}{z^n}\right).\ee
Clearly, the function $F=z^m$, where $\GCD(n,m)=1$, satisfies \eqref{homo} for some automorphism $\phi$, 
and corresponding diagram \eqref{dia3} has the form
$$T_m\circ \frac{1}{2}\left(z^n+\frac{1}{z^n}\right)= \frac{1}{2}\left(z^n+\frac{1}{z^n}\right)\circ z^m.$$
Further, $F$ maps the cyclic subgroup $\Gamma_{\f O^{\prime}}\subset \Gamma_{\f O}$, generated by $\beta$ and corresponding to the orbifold $\f O^{\prime}$ 
defined  by  the condition $\nu(-1)=\nu(1)=2$,
to itself, and  equality \eqref{iuy}
has the form 
$$\frac{1}{2}\left(z^n+\frac{1}{z^n}\right)=T_n\circ \frac{1}{2}\left(z+\frac{1}{z}\right).$$ Thus, we arrive 
to 
the diagram 
\begin{equation*} 
\begin{CD} 
\tt{\f O^{\prime}} @>z^m>> \tt{\f O^{\prime}}\\
@VV\frac{1}{2}\left(z+\frac{1}{z}\right)  V @VV \frac{1}{2}\left(z+\frac{1}{z}\right) V\\ 
\f O^{\prime} @> T_m >> \f O^{\prime}\\
@VV T_n  V @VV T_n  V\\ 
\f O @>T_m >> \f O\ 
\end{CD}
\end{equation*}
and to 
the series of semiconjugate functions 
\be \l{pila} T_m\circ T_n=T_n\circ T_m.\ee

\vskip 0.2cm
Examples of solutions of \eqref{i1} involving  Latt\`es functions can be constructed similarly. 
We start from a covering map $A:\f O\rightarrow \f O$ for some orbifold $\f O$ of zero Euler characteristic
and a function $F$ which makes the diagram 
\be \l{dia33}
\begin{CD}
\tt {\f O} @>F>> \tt {\f O}\\
@VV\theta_{\f O}V @VV\theta_{\f O}V\\ 
\f O @>A >> \f O\ 
\end{CD}
\ee
commutative (notice that since $A$ is a covering map, $F$ is an isomorphism and hence 
has the form $F=az+b,$ $a,b\in \C$). Since 
$\theta_{\f O}$ is transcendental, 
diagram \eqref{dia33} by itself does not provide now any rational solution of \eqref{i1}.  However, it is easy to see that if 
$\Gamma_{\f O\prime}\subseteq \Gamma_{\f O}$ is a subgroup corresponding to another orbifold $\f O^\prime$ with $\chi(\f O^{\prime})=0$ such that   
$\phi$ in \e{homm} satisfies the condition $\phi(\Gamma_{\f O\prime})\subseteq \Gamma_{\f O\prime}$, then there exists a rational function $X$ such that diagram  \e{gooopa2} is commutative and the corresponding Latt\`es functions $B$ and $A$  are semiconjugate.

For example, let  $\f O$ be an  orbifold with ramification 
$(2,2,2,2)$. For such an orbifold the group $\Gamma_{\f O}$ is generated by translations by elements of some lattice $L$ of rank two in $\C$ and the transformation $z\rightarrow -z$. The universal cover 
$\theta_{\f O}$ is the Weierstrass function $\wp_L$ corresponding to $L$.  Clearly, the function $F=mz,$ where $m\geq 2$ is an integer, satisfies \eqref{homm}.
The corresponding Latt\`es function $A=R_{L,m}$ is a rational function satisfying the equality  
$\wp_L(mz)=R_{L,m}\circ \wp_L$. 
Let now $L^{\prime}$ be 
a sublattice of $L$. Then $\wp_L=X\circ \wp_{L^{\prime}}$ for some 
rational function $X$, and $F(L^{\prime})\subset L^{\prime}$, implying that  
$$R_{L,m}\circ X=X\circ R_{L^{\prime},m}.$$

\vskip 0.2cm In conclusion, let us make several comments regarding the equivalence relation $\sim.$ 
First, observe that for a rational function $F$  its equivalence  class may simply coincide with its conjugacy class. 
This is the case for examples for any rational function $F$ which cannot be decomposed into a composition of rational functions of lesser degrees. Notice that in fact any generic rational function $F$  is indecomposable 
since the condition $F=A\circ B$ imposes algebraic restrictions on the coefficients of $F$.  
Another example of a function whose equivalence class coincides with its conjugacy class is the function $F=z^n$. Although this function is decomposable, any decomposition $A\circ B$ of $z^n$ has the form 
$$A=z^d\circ \mu, \ \ \ \ B=\mu^{-1}\circ z^{n/d}$$ for some $d\vert n$ and a M\"obius transformation $\mu,$ and therefore any transformation of the form $A\circ B\rightarrow B\circ A$ leads to a conjugated function.

Further, by the result of McMullen cited in the introduction, any equivalence class contains at most a finite number of conjugacy classes. However, there is no uniform bound 
on the number of such classes. 
Indeed, for any rational function $R$ and a natural number $d$ we have:
$$z^2\circ zR(z^{2^d})=zR^2(z^{2^{d-1}})\circ z^2,$$
$$z^2\circ zR^2(z^{2^{d-1}})=zR^4(z^{2^{d-2}})\circ z^2,$$
$$z^2\circ zR^4(z^{2^{d-2}})=zR^8(z^{2^{d-3}})\circ z^2,$$
$$...$$
$$z^2\circ zR^{2^{d-1}}(z^{2})=zR^{2^d}(z)\circ z^2,$$
implying that 
\be \l{bliii}  z^2\circ zR(z^{2^d})\sim z^2\circ zR^2(z^{2^{d-1}})\sim z^2\circ zR^4(z^{2^{d-2}})\sim \dots \sim  z^2\circ zR^{2^{d-1}}(z^2).\ee
Setting for example  $R=z-1$, it is easy to see that the corresponding  functions 
$$F_i=  z^2(z^{2^{d-i+1}}-1)^{2^i}, \ \ \ i=1,\dots d,$$
in \eqref{bliii}  
cannot be conjugated. 
Indeed, since $F_i^{-1}\{\infty\}=\{\infty\}$, $1\leq i \leq d,$
the equality 
\be \l{meba} F_i=\mu \circ F_j\circ \mu^{-1},\ \ \ i\neq j,\ee where $\mu$ is a M\"obius transformation, implies  that 
$\mu=\alpha z+\beta,$ $\alpha,\beta\in \C.$ Furthermore, denoting by $n$  the common degree of the functions $F_i$, $1\leq i \leq d,$ and 
comparing the coefficients of $z^{n-1}$ 
in both parts of equality 
\eqref{meba}, we conclude that $\beta=0.$ 
Finally, \eqref{meba} cannot be satisfied for $\mu=\alpha z$ since polynomials in the left and in the right parts of \eqref{meba} have 
different collections of monomials with non-zero coefficients.

\vskip 0.2cm
\noindent{\bf Acknowledgments}. The author is grateful to 
A. Eremenko and C. McMullen for discussions, and to the Max-Planck-Institut fuer Mathematik for the hospitality and the support. Besides, the author would like to thank the referee for many valuable 
comments.
% and for pointing out the link between the equivalence $\sim$ and isospectral rational functions.

\vskip 0.4cm
\centerline{ADDENDUM}

\vskip 0.2cm
The very first version of this paper was published in the form of a preprint in 2011. Below
we mention some publications related to  semiconjugate rational functions that appeared afterward.

In the paper \cite{e} functional  equation \eqref{i1} was related with the problem of describing of Jordan curves in $\C$ invariant under a rational function.
On the other hand, it was shown in the  paper \cite{ms} that equation \eqref{i1} is closely related to the  problem of describing of algebraic
curves $\f C$ in $\C^2$ invariant under maps of the form $F:\,(x,y)\rightarrow (f(x),g(y)),$ where $f,g$ 
are polynomials of degree at least two. In particular, the paper   \cite{ms} contains a detailed analysis of equation \eqref{i1} in the polynomial case, 
based on the Ritt theory of decompositions of polynomials.

Another approach to equation \eqref{i1} in the polynomial case, using results of \cite{p1}
about polynomials sharing preimages of compact sets, was proposed in  \cite{pj}. In particular, it was shown in  \cite{pj} that in the polynomial case 
the conditions $A\leq B$ and $B\leq A$ hold simultaneously if and only if 
$A\sim B$. Notice however that methods of both papers \cite{ms} and \cite{pj} seem to be restricted to the polynomial case.

It the paper \cite{gen} the methods of this paper were applied to the functional equation 
$ A\circ C=D\circ B,$ 
where $A,B,C,D$ are rational functions.
In particular, it was shown in \cite{gen} that if  $A,B,C,D$ satisfy this equation and 
the algebraic curve $$A(x)-D(y)=0$$ is irreducible, then whenever $\deg D\geq 84\, \deg A$
the inequality $\chi(\f O_2^A)\geq 0$ holds.  
Finally, further development of ideas and methods of this paper was given in the paper \cite{pnew} devoted to  quantitative aspects of the description of 
solutions of \eqref{i1} for fixed $B$. In particular, it was shown in \cite{pnew} that  if $B$ is neither a Latt\`es function nor conjugated to $z^{\pm d}$ or $\pm T_d$, 
then, up to some natural transformations, the number of $A$ and $X$ satisfying \eqref{i1} is finite  and can be effectively bounded in terms of $\deg B$ only.

\end{document}